\documentclass[11pt]{article}
\input epsf
\usepackage{amsfonts}
\usepackage{amscd}
\usepackage{amsmath,amssymb}
\usepackage{amstext}
\usepackage{amscd}
\usepackage[all]{xy}
\usepackage{pstricks,pst-node,pst-tree}

\newtheorem{lem}{LEMMA}[section]
\newtheorem{theo}[lem]{THEOREM}

\newtheorem{coro}[lem]{COROLLARY}
\newtheorem{prop}[lem]{PROPOSITION}

\newtheorem{definition}[lem]{Definition}

\newtheorem{rem}[lem]{Remark}

\newtheorem{exs}[lem]{Examples}
\renewcommand{\descriptionlabel}[1]%
       {\hspace{\labelsep}\textsf{#1}}

\begin{document}

\title{Wild knots in higher dimensions as limit sets of  Kleinian groups}
\author{Margareta Boege\thanks{Research partially supported by PFAMU-DGAPA.}, Gabriela Hinojosa\thanks{Research  partially supported by CONACyT CB-2007/83885.},
Alberto Verjovsky\thanks{Research partially
supported by CONACyT proyecto U1 55084, and PAPIIT (Universidad
Nacional Aut\'onoma de M\'exico) \# IN102108.}}
\maketitle

\begin{abstract}In this paper we construct infinitely many wild knots,
$\mathbb{S}^{n}\hookrightarrow\mathbb{S}^{n+2}$, for $n=1,2,3,4$ and $5$, each of which
is a limit set of a geometrically finite Kleinian group. We also describe some of their properties.
\end{abstract}

\section{Introduction}

 Kleinian groups were introduced by Henri Poincar\'e  in the 1880's  \cite {Po},  as
  the monodromy groups of certain second order differential
  equations on the Riemann sphere $\widehat{\mathbb{C}}$. They
  have  played a major  role in many parts of mathematics throughout
  the twentieth and the present centuries,
as for example in Riemann surfaces and Tei\-chm\"uller theory, automorphic
  forms,
holomorphic dynamics, conformal and hyperbolic geometry, number theory, and to\-pology
(for instance the study of 3-manifolds).\\

The higher dimensional analogue of Kleinian groups are certain discrete
subgroups of the group of diffeomorphisms of the $(n+2)$-sphere ($n\geq1$), with its standard metric,
consisting of those diffeomorphisms which
preserve angles and which we denote by ${M\ddot{o}b}\,(\mathbb S^{n+2})$. The subgroup of index two of
${M\ddot{o}b}\,(\mathbb S^{n+2})$
which consists of those elements which are orientation-preserving is
called the {\it Conformal or M\"obius Group } and is denoted by
${M\ddot{o}b}_+\,(\mathbb S^{n+2})$. See the book by Ahlfors \cite{ahlfors}. If
$\Gamma\subset {M\ddot{o}b}\,(\mathbb S^{n+2})$ is a discrete subgroup
acting conformally on the $(n+2)$-sphere then this action extends naturally to a conformal action on the disk $\mathbb{D}^{n+3}$. Its limit set, $\Lambda(\Gamma)$,
is the set of points of $\mathbb{S}^{n+2}$ which are accumulation points of some orbit of
$\Gamma$ in $\mathbb{D}^{n+3}$. If
$\Omega(\Gamma):=\mathbb{S}^{n+2}-\Lambda(\Gamma)\neq\emptyset$
one says that $\Gamma$ is a {\it Kleinian group}. The set $\Omega(\Gamma)$ is called the {\it discontinuity set} of
$\Gamma$.

One interesting question
is whether a topological $n$-sphere ($n\geq1$) which is not a round sphere can be the limit set of a higher dimensional
{\it geometrically finite} Kleinian group. In
this case one can show that the sphere is necessarily
fractal (possibly unknotted). Examples of wild knots in $\mathbb S^3$,
which are limit sets of geometrically finite  Kleinian groups, were obtained by
 Maskit (\cite{maskit}), Kapovich (\cite{kap1}), Hinojosa
(\cite{hinojosa}) and Gromov, Lawson and Thurston (\cite{GLT}).
An example of a wild 2-sphere in $\mathbb{S}^4$ which is the limit set of a geometrically
finite Kleinian group was obtained by
 the second-named author \cite{hinojosa2} and, independently, by Belegradek \cite{belegradek} (see also \cite{apanasov} for a wild 
limit set $\mathbb{S}^{2}\rightarrow \mathbb{S}^{3}$) .
Such wild knots are examples of self-similar fractal sets and they are extremely beautiful
to contemplate. For instance, one can admire the pictures in the classic book by R. Fricke and F. Klein (\cite{klein})
or the pictures of limit sets of classical Kleinian groups
in the book {\it Indra's Pearls: The vision of Felix Klein} by  D. Mumford, C. Series, D. Wright (\cite{mumford}).

In this paper we construct an infinite number of wild knots $\mathbb S^{n}\hookrightarrow\mathbb S^{n+2}$
for $n=1,\cdots ,5$
which are limit sets of geometrically finite Kleinian groups.
If there existed a way to picture and travel through the spheres of high dimensions we could contemplate our examples as
the analogue of Indra's Pearls in higher dimensions!

This paper is organized as follows. In section 2 we give the preliminaries: the definition of oriented tangles
and knots in high dimensions and some basic facts of Kleinian groups.
In section 3 we introduce the notion of {\it orthogonal ball covering} (OBC) of $\mathbb{R}^{n+2}$ as a covering by
round balls satisfying certain conditions of orthogonality. Using results by  L. Potyagailo, E.B. Vinberg in \cite{PV}, 
we describe explicit OBCs for $n=1,\ldots,5$. In section 4
we construct $n$-knots as limit sets of geometrically finite Kleinian groups for $n=1,\ldots,5$. In section 5 we show that if we start with a
nontrivial tame fibered $n$-knot $K$, then the complement of the
corresponding limit $n$-knot $\Lambda(K)$ also fibers over $\mathbb{S}^1$ and give a description of the fibers.
In section 6 we describe the monodromy of  $\Lambda(K)$ in terms of the monodromy of $K$. We prove that  $\Lambda(K)$ is wildly embedded in
$\mathbb{R}^{n+2}$. Therefore there exist infinitely many wild knots which are limit sets of discrete,
 geometrically finite Kleinian groups in dimensions 3, 4, 5, 6 and 7.

\section{Preliminaries}

In classical knot theory, a subset $K$ of a space $X$ is a {\it knot} if $K$ is homeomorphic to a sphere
$\mathbb{S}^{p}$. Two knots $K$, $K'$ are {\it equivalent} if there is a homeomorphism $h:X\rightarrow X$ such that $h(K)=K'$;
in other words $(X,K)\cong (X,K')$. However, a knot $K$ is sometimes defined to be an embedding
$K:\mathbb{S}^{p}\rightarrow\mathbb{S}^{n}$ (see \cite{mazur}, \cite{rolfsen}).
We shall also find this convenient at times and will use the same symbol to
denote either the map $K$ or its image $K(\mathbb{S}^{p})$ in $\mathbb{S}^{n}$.

\begin{definition}
An {\it oriented n-dimensional tame single-strand tangle} is a couple $D=({B}^{n+2},T)$ satisfying the following conditions:
\begin{enumerate}
\item ${B}^{n+2}$ is homeomorphic to the (n+2)-disk $\mathbb{D}^{n+2}$, and $T$ is homeomorphic to the n-disk $\mathbb{D}^{n}$.
\item The pair $({B}^{n+2},T)$ is a proper manifold pair, i.e. $\partial T\subset \partial {B}^{n+2}$ and
$\mbox{Int}(T)\subset \mbox{Int}({B}^{n+2})$.
\item $(B^{n+2},T)$ is locally flat (\cite{rushing}, p.33).
\item ${B}^{n+2}$ has an orientation which induces the canonical orientation on its boundary
$\partial B^{n+2}$.
\item $(\partial B^{n+2},\partial T)$ is homeomorphic to $(\partial\mathbb{D}^{n+2},\partial\mathbb{D}^{n})$.
\end{enumerate}
\end{definition}

Compare to Zeeman's definition of ball-pair in (\cite{zeeman}).\\

Two oriented tangles $D_{1}=({B}_{1}^{n+2},T_{1})$, $D_{2}=({B}_{2}^{n+2},T_{2})$ are equivalent if there exists an orientation-preserving
homeomorphism of ${B}_{1}^{n+2}$ onto ${B}_{2}^{n+2}$ that sends $T_{1}$ to $T_{2}$. A tangle is unknotted if it is
equivalent to the trivial tangle $(\mathbb{D}^{n+2},\mathbb{D}^{n})$.\\

Given an oriented tangle $D=({B}^{n+2},T)$, the pair $(\partial B^{n+2},\partial T)$ is homeomorphic to the pair
$(\mathbb{S}^{n+1},\mathbb{S}^{n-1})$, via a homeomorphism $f$.
Then $D$ determines canonically a knot $K\subset\mathbb{S}^{n+2}$, in the following way:
$(\mathbb{S}^{n+2},K)=({B}^{n+2},T)\cup_{f}(\mathbb{D}^{n+2},\mathbb{D}^{n})$.\\

Conversely, given a smooth knot $K\subset\mathbb{S}^{n+2}$, there exists a smooth ball $B^{n+2}$ such that $(B^{n+2},B^{n+2}\cap K)$
is equivalent to the trivial tangle. The tangle
$K_{T}=(\mathbb{S}^{n+2}-\mbox{Int}(B^{n+2}),K-\mbox{Int}(B^{n+2}\cap K))$ is called the {\it
canonical tangle} associated to $K$. Notice that if $K$ is not the trivial knot, then $K_{T}$ is not equivalent to
the trivial tangle. In this case, we say that $K_{T}$ is knotted.

\begin{rem}
The above constructions are  well-defined up to isotopy.
\end{rem}

The {\it connected sum} of the oriented tangles $D_{1}=({B}_{1}^{n+2},T_{1})$ and $D_{2}=({B}_{2}^{n+2},T_{2})$ for $n>2$,
denoted by $D_{1}\# D_{2}$,
can be defined as follows:  Since  $D_{i}$ is locally flat $i=1,2$, there exist sets $U_{i}\subset \partial B_{i}$ closed in $\partial B_{i}$
such that $\mbox{Int}(U_{i})\cap \partial T_{i}\neq\emptyset$ and the pair $(U_{i}, U_{i}\cap T_{i})$ is homeomorphic to
$(\mathbb{D}^{n+1},\mathbb{D}^{n-1})$. Choose an orientation-reversing homeomorphism $h$ of
$(U_{1}, U_{1}\cap T_{1})$ onto
$(U_{2}, U_{2}\cap T_{2})$. Then
$$
D_{1}\# D_{2}=({B}_{1}^{n+2},T_{1})\cup_{h} ({B}_{2}^{n+2},T_{2})
$$
\begin{rem}
The  connected sum does not depend on the choice of the homeomorphism  $h$ and the open sets $U_{i}$.
\end{rem}

Our goal is to obtain wild $n$-spheres, $n=1,2,\ldots, 5$ as limit sets of conformal
Kleinian groups. We will briefly review some basic definitions about Kleinian
groups.\\

Let $M\ddot{o}b\,(\mathbb{S}^{n})$ denote the group of M\"obius
transformations of the $n$-sphere
$\mathbb{S}^{n}=\mathbb{R}^{n}\cup\{\infty\}$, {\it i.e.}, conformal diffeomorphisms of
$\mathbb{S}^{n}$ with respect to the standard metric. For a discrete group
$G\subset M\ddot{o}b(\mathbb{S}^{n})$ the {\it discontinuity set} $\Omega(G)$ is
defined as follows:
$$
\Omega (G)=\{x\in\mathbb{S}^{n}: \mbox{the point}\hspace{.2cm}
x\hspace{.2cm} \mbox{possesses a neighborhood}
\hspace {.2cm}U(x)\hspace{.2cm}\mbox{such that}
$$
$$
\hspace{.6cm}
U(x)\cap g(U(x))\hspace{.2cm}
\mbox{is empty for all but a finite number of elements}\hspace{.2cm} g\in G\}.
$$
The complement
$\mathbb{S}^{n}-\Omega(G)=\Lambda(G)$ is called the {\it  limit set} (see
\cite{kap1}). Both $\Omega(G)$ and $\Lambda(G)$ are $G$-invariant. The set $\Omega(G)$ is
open, hence $\Lambda(G)$ is compact.\\

A subgroup $G\subset M\ddot{o}b\,(\mathbb{S}^{n})$ is called {\it  Kleinian} if
$\Omega (G)$ is not empty.\\

We recall that a conformal map $\psi$ on $\mathbb{S}^{n}$ has a  Poincar\'e extension
to  the hyperbolic space $\mathbb{H}^{n+1}$ as an isometry with
respect to the Poincar\'{e} metric. Hence we can identify
the group $M\ddot{o}b(\mathbb{S}^{n})$ with the group of
isometries of hyperbolic $(n+1)$-space $\mathbb{H}^{n+1}$. This
allows us to define the limit set of a Kleinian group through sequences. We say that
a point $x$ is a {\it  limit point} for the Kleinian group
$G$, if there exists a point $z\in\mathbb{S}^{n}$ and a sequence $\{g_{m}\}$
of {\it distinct elements} of $G$, with $g_{m}(z)\rightarrow x$.
The set of limit points is $\Lambda (G)$ (see \cite{maskit} section II.D).\\

One way to illustrate the action of a Kleinian group $G$ is to draw a picture
of $\Omega(G)/G$. For this purpose a fundamental domain is very helpful. Roughly
speaking, it contains one point from each equivalence class in $\Omega(G)$ (see
\cite{kap2} pages  78-79, \cite{maskit} pages  29-30).\\

\begin{definition}
A fundamental domain $D$ for a Kleinian group $G$ is a
co\-di\-men\-sion zero piecewise-smooth submanifold (sub polyhedron) of
$\Omega(G)$ satisfying the following
\begin{enumerate}
\item $\bigcup_{g\in G} g(\mbox{Cl}_{\Omega(G)} D)=\Omega$ where $Cl_{\Omega(G)}$ is the closure in $\Omega(G)$.
\item $g(\mbox{Int}(D))\cap \mbox{Int}(D)=\emptyset$ for all
$g\in G-\{e\}$, where $e$ is the identity in $G$ and $\mbox{Int}$ denotes the interior.
\item The boundary of $D$ in $\Omega (G)$ is a piecewise-smooth
(polyhedron) submanifold in $\Omega (G)$, divided into a union
of smooth submanifolds (convex polyhedra) which are called {\rm
faces}. For
each face $S$, there is a corresponding face $F$ and an element $g_{SF}\in
G-\{e\}$ such that $g_{SF}(S)=F$ ($g_{SF}$ is called a  face-pairing
transformation); $g_{SF}=g_{FS}^{-1}$.
\item Only finitely many translates of $D$ meet any compact subset of
$\Omega (G)$.
 \end{enumerate}
\end{definition}

\begin{theo}(\cite{kap2}, \cite{maskit}) Let
$D^{*}=\overline{D}\cap\Omega/\sim_{G}$ denote the orbit space with
the quotient topology. Then $D^{*}$ is homeomorphic to $\Omega/G$.
\end{theo}

\section{Orthogonal ball coverings of $\mathbb{R}^{n+2}$}

A countable collection of closed round $(n+2)$-balls $B_{1}$, $B_{2}$,
$B_{3},\ldots$ is an {\it orthogonal ball covering} (OBC) of $\mathbb{R}^{n+2}$ if
\begin{enumerate}
\item $\cup_{i} B_{i}=\mathbb{R}^{n+2}$
\item There exist  $\epsilon_{1}, \epsilon_{2}>0$ such that  $\epsilon_{1}<\mbox{diameter} (B_{i})\leq \epsilon_{2}$, $\forall i$
\item For each pair $(B_{i},B_{j})$ with $i\neq j$, we have that they are either disjoint, meet at only one point
or their boundaries meet orthogonally.
\end{enumerate}

In particular, an OBC is locally finite.\\

One has the following theorem of Potyagailo and Vinberg (see \cite{PV}):

\begin{theo}
There exist right-angled hyperbolic polyhedra of finite volume, with at least one point at infinity, in
$\mathbb{H}^{n}$,  for $n=3,\ldots,8$.
\end{theo}

This theorem has, for our purposes, the following relevant corollary:

\begin{coro}
There exist OBCs for $\mathbb{R}^{n}$, $n=2,\ldots,7$.
\end{coro}

{\it Proof.}
The OBC is obtained as follows: Consider the
half-space model of $\mathbb{H}^{n+1}=\{(x_1,\cdots,x_{n+1})|\,\,\, x_{n+1}>0\}\subset \mathbb R^{n+1}$ with boundary the hyperplane
at infinity $\mathbb R^{n}$ with equation $x_{n+1}=0$. Let $P$ be a right-angled hyperbolic polyhedra of finite volume in
 $\mathbb{H}^{n+1}$ for  $n=2,\ldots,7$ as in the previous theorem.
 We can assume that one
of the ideal vertices of $P$ is $\infty$. Consider the facets $F_{i}$ of $P$ which are not asymptotic to $\infty$ and let the spheres
$S_{i}$ ($i=1,\ldots, k$) be the ideal boundaries of the hyperbolic hyperplanes through $F_{i}$'s. Now, apply to the spheres $S_{i}$ the
group $G$ generated by the reflections in the facets of $P$ asymptotic to $\infty$ (semi-hyperplanes orthogonal to the hyperplane at
infinity). The action of $G$ on the hyperplane at infinity has as fundamental domain a compact parallelepiped (in fact, as we will see 
afterwards, it is a regular cube). The set of images of
the balls whose boundaries are the spheres $S_{i}$ is the desired OBC.

\subsection{Description of some OBCs for $\mathbb{R}^{n}$, $n=2,\ldots,7$ and their nerves}

We will describe geometrically the OBCs $\cal {B}$ obtained in the corollary 3.2 and the geometric realization of their nerves (compare M. Kapovich \cite{kap3}).
Given a collection of round open balls $\{B_{j}, j\in I\}$ in
$\mathbb{R}^{n}$, consider its nerve $N$. We define the canonical simplicial mapping $f:N\rightarrow \mathbb{R}^{n}$
by sending each vertex of $N$ to the center of the corresponding ball and extending $f$ linearly to the simplices of
$N$. {\it The geometric realization} of  $N$ is $f(N)$.\\

For $n=2$, we have that the fundamental domain for the group $G$ is the unit square, $I^{2}$. We set two closed round balls of
radius one at $(0,0)$ and $(1,1)$ (see Figure 1). The geometric realization of its nerve, is the closed segment joining
$(0,0)$ and $(1,1)$.\\

\begin{figure}[tbh]
\centerline{\epsfxsize=1.6in \epsfbox{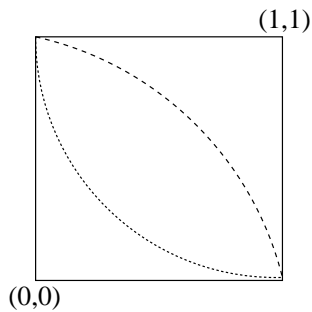}}
\caption{\sl A fundamental domain for $n=2$.}
\label{F1}
\end{figure}

If we propagate the fundamental domain by $G$, we get a flower at the remaining vertices of $I^{2}$, {\it  i.e.,} we obtain a
configuration consisting of four balls such that balls centered at the same straight line parallel to a coordinate axis, are
tangent (see Figure 2). The geometric realization of the nerve of the above four balls is a rhombus.\\

\begin{figure}[tbh]
\centerline{\epsfxsize=1.6in \epsfbox{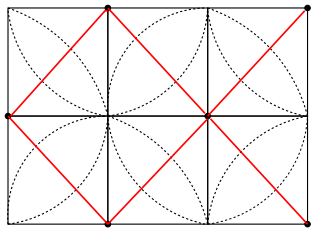}}
\caption{\sl A flower for $n=2$.}
\label{F2}
\end{figure}

For $n=3$, we have that the fundamental domain for the group $G$ is the unit cube, $I^{3}$. Notice that $I^{3}=I^{2}\times [0,1]$.
At the face $I^{2}\times \{0\}$, we set two closed round balls of radius one at $(0,0,0)$ and $(1,1,0)$. At the face
$I^{2}\times \{1\}$, we set two closed round balls of radius one at $(1,0,1)$ and $(0,1,1)$ (see Figure 3). The
geometric realization of its nerve is a solid tetrahedron spanned by $(0,0,0)$, $(1,1,0)$, $(1,0,1)$ and $(0,1,1)$.\\

\begin{figure}[tbh]
\centerline{\epsfxsize=1.6in \epsfbox{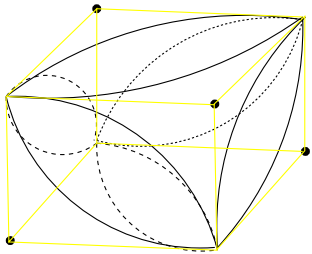}}
\caption{\sl A fundamental domain for $n=3$.}
\label{F3}
\end{figure}

As in the previous case, if we propagate the fundamental domain by $G$, we get a flower at the remaining vertices of $I^{3}$,
{\it  i.e.,} we obtain a
configuration consisting of six balls such that balls centered at the same straight line parallel to a coordinate axis are
tangent (see Figure 4). This situation will appear in the next dimensions, therefore we will give a general definition

\begin{figure}[tbh]
\centerline{\epsfxsize=1.6in \epsfbox{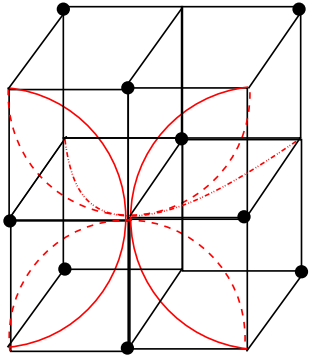}}
\caption{\sl A flower for $n=3$.}
\label{F4}
\end{figure}

\begin{definition}
Let $B_{1},B_{2},\ldots,B_{2n+4}$ be round closed balls in the OBC $\cal{B}$ of $\mathbb{R}^{n+2}$. We say that they form a flower if
\begin{enumerate}
\item For each pair $(B_{i},B_{j})$ with $i\neq j$, we
have that they either meet at only one point or their boundaries meet orthogonally.
\item The intersection $\cap_{i=1}^{2n+4} B_{i}$ consists of only one point $c$ which we call the {\it center} of the flower.
\end{enumerate}
\end{definition}

The geometric realization of the nerve of the above six balls is the boundary of a solid octahedron. Observe
that each face of the octahedron corresponds to a face of a solid tetrahedron previously described.\\

For $n=4$, we have that the fundamental domain for the group $G$ is the unit hypercube, $I^{4}$. Notice that $I^{4}=I^{3}\times [0,1]$.
At the face $I^{3}\times \{0\}$, we set four closed round balls of radius one at $(0,0,0,0)$, $(1,1,0,0)$, $(1,0,1,0)$ and $(0,1,1,0)$.
At the face
$I^{3}\times \{1\}$, we set four closed round balls of radius one at $(1,1,1,1)$, $(0,1,0,1)$, $(1,0,0,1)$ and $(0,0,1,1)$. \\

Observe that all these balls meet at the center of  $I^{4}$, $(\frac{1}{2},\frac{1}{2},\frac{1}{2},\frac{1}{2})$,
hence they form a flower. The geometric realization of the nerve of the above eight balls is the 4-dimensional analogous of
the octahedron, a hyper-octahedron, which for simplicity will be called a {\it diamond}. Notice that the center of
this diamond is the center of the corresponding flower. This remains true in the next dimensions.\\

For $n=5$  we have that the fundamental domain for the group $G$ is the unit hypercube, $I^{5}$. Notice that $I^{5}=I^{4}\times [0,1]$.
In the face $I^{4}\times \{0\}$, we set eight closed round balls of radius one centered at the following vertices
$(0,0,0,0,0)$, $(1,1,0,0,0)$, $(1,0,1,0,0)$,
$(0,1,1,0,0)$, $(1,1,1,1,0)$, $(0,1,0,1,0)$, $(1,0,0,1,0)$ and $(0,0,1,1,0)$.
In the face
$I^{4}\times \{1\}$, we set eight closed round balls of radius one at
$(1,0,0,0,1)$, $(0,1,0,0,1)$, $(0,0,1,0,1)$, $(0,0,0,1,1)$, $(0,1,1,1,1)$, $(1,0,1,1,1)$, $(1,1,0,1,1)$ and $(1,1,1,0,1)$. \\

Next, we set a ball at $(\frac{1}{2},\frac{1}{2},\frac{1}{2},\frac{1}{2},\frac{1}{2})$ of radius $\frac{1}{2}$.\\

If we propagate the fundamental domain by $G$, as before we obtain flowers at vertices and also a flower at the center of each
hyper-face of dimension four.\\

For $n=6$  we have again that the fundamental domain for the group $G$ is the unit hypercube, $I^{6}$. At each face of $I^{6}$, we
 repeat the construction for  $n=5$ in such a way that, we set
32 closed round balls of radius one centered at vertices and 12 closed round balls of radius $\frac{1}{2}$, each one set at the
center of each face.\\

Notice that the twelve closed balls at the centers of faces, form a flower at the center of $I^{6}$. \\

For $n=7$ we have again that the fundamental domain for the group $G$ is the unit hypercube, $I^{7}$. We will not describe in detail this
construction, since we only need for our purpose, the OBC restricted to its faces. Now,  at each face of $I^{7}$, we
 repeat the construction for  $n=6$. \\

\begin{rem}
The geometric realization of the nerve of the OBC $\cal{B}$ is embedded in $\mathbb{R}^{n+2}$, for $n=1,\ldots,5$.
\end{rem}

\section{Construction of a wild knot as limit set of a Kleinian Group}

As we have mentioned before, from the existence of  right-angled hyperbolic polyhedra of finite volume, with at least
one point at infinity in $\mathbb{H}^{n+3}$,  for $n=1,\ldots,5$, follows that there exist OBC, ${\cal{B}}$, for $\mathbb{R}^{n+2}$,
$n=1,\ldots,5$. In this section, we will use these OBCs to construct  wild knots. Our aim, in this and the next sections, is to prove the following:

\begin{theo}
There exist infinitely many non-equivalent knots $\psi:\mathbb{S}^{n}\rightarrow \mathbb{S}^{n+2}$ wildly embedded
as limit sets of geometrically finite Kleinian groups, for $n=1,2,3,4,5$.
\end{theo}

Let $\psi:\mathbb{S}^{n}\rightarrow \mathbb{R}^{n+2}\subset \mathbb{R}^{n+2}\cup\{\infty\}=\mathbb{S}^{n+2}$
be a smoothly embedded knotted $n$-sphere in $\mathbb{S}^{n+2}$. We denote $K=\psi(\mathbb{S}^{n})$ and endow it with the
Riemannian metric induced by the standard Riemannian
metric of $\mathbb{R}^{n+2}$.\\

The general idea of our construction of wild knots is the following: Given a smooth knot $\bar{K}$, there
exists an isotopic copy of it, $K$, in the $n$-skeleton of the canonical cubulation of $\mathbb{R}^{n+2}$
 (see the next subsection). This cubulation is canonically associated to an OBC, ${\cal{B}}$.

Let ${\cal{B}}(K)=\{B\in {\cal{B}}\,:\,K\cap\mbox{Int}(B)\neq\emptyset\}$. The group
generated by inversion on the boundaries of balls belonging to ${\cal{B}}(K)$ whose centers are in the $n$-skeleton will be Kleinian and its limit set will be a wild knot.\\

\subsection{OBCs and Cubulations for $\mathbb{R}^{n+2}$}

A {\it cubulation} of $\mathbb{R}^{n+2}$ is a decomposition of $\mathbb{R}^{n+2}$ into a collection $\cal C$ of $(n+2)$-dimensional
cubes such that any two of its hypercubes are
either disjoint or meet in one common face of some dimension. This provides  $\mathbb{R}^{n+2}$ with the structure of a cubic
complex.\\

The canonical cubulation $\cal C$ of $\mathbb{R}^{n+2}$ is the decomposition into hypercubes which are the images of the unit cube
$I^{n+2}=\{(x_{1},\ldots,x_{n+2})\,|\,0\leq x_{i}\leq 1\}$ by translations by vectors with integer coefficients.
Consider the homothetic transformation $H_{m}:\mathbb{R}^{n+2}\rightarrow\mathbb{R}^{n+2}$,
$H_m(x)=\frac{1}{m} x$,
where $m>1$ is an integer. The set ${\cal{C}}_{m}=H_{m}(\cal{C})$ is called a {\it subcubulation} or {\it cubical subdivision} of $\cal C$.\\

Observe that the $n$-skeleton of $\cal C$, denoted by $\cal S$, consists of the union of the $n$-skeletons of the cubes in  $\cal C$,
{\it i.e.,} the union of all cubes of dimension $n$ contained in the faces of the $(n+2)$-cubes in  $\cal C$.
We will call $\cal S$ the {\it standard scaffolding} of $\mathbb{R}^{n+2}$.\\

 By the previous section, the OBCs of ${\Bbb R}^{n+2}$
were obtained by covering the unit cube by balls which are either tangent, meet orthogonally or are disjoint and then
taking all the balls obtained  as images under the group generated by reflections on the
faces of the cube. This is possible since the balls meet the planes which support the faces either tangentially or orthogonally.
Thus we have the canonical cubulation associated in a natural way to the OBCs.\\

We can also associate the cubulation ${\cal C}_{2}$ to the OBCs. In this case, our fundamental cube $E$ is the union of $2^{n+2}$ cubes in $\cal C$ and
can be obtained as follows: take a cube $I$ in the standard cubulation and take a vertex $v\in I$ and reflect in all the hyperplanes which support
faces not containing $v$. The union of all the images of $I$ under these reflections is $E$. \\

Notice that there are two types of vertices in $I$: Those which are centers of balls and those which are not. We will choose $v$ such that it is not a center of
ball. Since the vertices of $E$ are the orbit of $v$, then all vertices of $E$ are not centers of balls. We can assume (if necessary after applying a translation) that
the origin is a vertex of $E$.\\

We will use the following theorem \cite{BHV} to embed an isotopic copy of the $n$-knot $K$ in the scaffolding ${\cal{S}}_{2}$ of ${\cal C}_{2}$.

\begin{theo} Let $\cal{C}$ be the standard cubulation of $\mathbb{R}^{n+2}$. Let $K\subset \mathbb{R}^{n+2}$ be a smooth knot of dimension $n$.
There exists a knot
$\hat{K}$ isotopic to $K$, which is contained in the scaffolding ($n$-skeleton) of the standard cubulation  $\cal{C}$
of $\mathbb{R}^{n+2}$. The cubulation of the knot $\hat{K}$ admits a subdivision by simplexes and
with this structure the knot is PL-equivalent to the $n$-sphere with its canonical PL-structure.
\end{theo}

By abuse of notation, we will denote $\hat{K}$  by $K$.

\subsection{Nerves and regular neighborhoods}

Let us consider the OBC $\cal {B}$ for $\mathbb{R}^{n+2}$ (section 3) together with the cubulation ${\cal C}_{2}$ constructed above. By theorem
4.2 we will assume that the knot $K$ is embedded in the scaffolding ${\cal{S}}_{2}$ of ${\cal C}_{2}$.\\

Let $o$ be a vertex of ${\cal C}_{2}$. By the above, we can assume that $o$ is the origin.
For $i,j\in\,\{1,2,,\ldots,n+2\}$, let $P_{i,j}=\{(x_{1},\ldots,x_{n+2})\,:\,x_{i}=x_{j}=0\}$
be the $n$-dimensional plane.

\begin{prop}
Let $B_{j}$ and  $B_{k}$ be balls whose centers lie in $P_{i,j}$ and $P_{i,k}$, respectively. Suppose that $B_{j}$ and  $B_{k}$ do not belong to a flower whose
center lies in $P_{i,j}\cap P_{i,k}$ (see definition 3.3). Then $B_{j}\cap B_{k}=\emptyset$.
\end{prop}

{\it Proof.}  Observe that  $P_{i,j}\cap P_{i,k}=\{(x_{1},\ldots,x_{n+2})\,:\,x_{i}=x_{j}=x_{k}=0\}$.
Let $C_{j}=(c^{j}_{1},c^{j}_{2},\ldots,c^{j}_{n+2})$,  $C_{k}=(c^{k}_{1},c^{k}_{2},\ldots,c^{k}_{n+2})$
be the centers of   $B_{j}$ and  $B_{k}$ respectively,  then $c^{j}_{i}=c^{j}_{j}=0$ and $c^{k}_{i}=c^{k}_{k}=0$. The minimum of the distance between
their centers is attained as $c^{j}_{r}=c^{k}_{r}$ for $r\neq i,j,k$. It is enough to prove this proposition for this case. \\

Let
$x=(x_{1},\ldots,x_{n+2})\in P_{i,j}\cap P_{i,k}$ such that $x_{r}=0$ if $r=i,j,k$  and $x_{r}=c^{j}_{r}=c^{k}_{r}$ if $r\neq i,j,k$. For $n<5$ the centers
$C_{j}$ and $C_{k}$ are vertices of the cubulation  ${\cal C}$ therefore their coordinates are integers.  Notice that $c^{j}_{k}$ and $c^{k}_{j}$
are bigger than one, since in other case
 $B_{j}$ and  $B_{k}$ would be part of the flower centered at  $x$. Hence these coordinates are either bigger or equal to two. This implies that the distance
$d(C_{j},C_{k})\geq\sqrt{8}$. Therefore $B_{j}\cap B_{k}=\emptyset$ since their radii are one.\\

For $n=5$ we have two types of balls: balls centered at vertices of the cubulation ${\cal C}$ and radii equal to one and,
balls whose centers coincide with centers of five dimensional faces of cubes in $\cal{C}$ and their radii are equal to $\frac{1}{2}$. The argument
to prove this case is analogous to the previous one. $\blacksquare$\\

In the remaining of this and the next sections we will restrict ourselves to those balls in  ${\cal{B}}$ which cover $K$. These balls coincide with the
balls whose centers lie in the scaffolding ${\cal{S}}_{2}$ of ${\cal C}_{2}$. That is,  let ${\cal{B}}(K)=\{B\in{\cal{B}}\,:\,K\cap \mbox{Int}(B)\neq\emptyset\}$
and ${T}=\cup\{B\in {\cal{B}}(K)\}$. Observe that $T$ is optimal in the sense that if we remove one, then the remaining balls do not cover $K$.
{\it We will call such a $T$ a generalized pearl necklace.} \\

Let $E$ be a cube in ${\cal{C}}_{2}$. Suppose that $K\cap E$ consists of more than one $n$-dimensional face  of $E$. Let
$F_{1}$ and $F_{2}$ be $n$-faces of $E$ such that $F_{1},\, F_{2}\subset K\cap E$. If $F_{12}$ denotes the $(n-1)$-dimensional cubic simplex
$F_{1}\cap F_{2}$, then {\it we will say that $K$ turns in $F_{12}$}.\\

From the construction of the OBCs we have that there are centers of flowers, $C_{i}$, $i=1,\ldots,r$, contained in $K$. At each center $C_{i}$
which does not belong to an $(n-1)$-cubic complex where $K$ turns, we will add a ball $B_{C_{i}}$ of center $C_{i}$ and radius less than $\frac{1}{\sqrt{n+2}}$.\\

Let ${\cal{B}}'(K)=\{B_{C_{i}}\}^{r}_{i=1}\cup {\cal{B}}(K)$ and $\tilde{T}=\cup_{i=1}^{r} B_{C_{i}}\cup T$.
Observe that the pair $(B,B\cap K)$, $B\in {\cal{B}}'(K)$
is isotopic to the trivial tangle. {\it We will call a $\tilde{T}$ a generalized increased pearl necklace.}\\

Let ${\cal{B}}_{K}=\{K\cap \mbox{Int}(B_{i})|B_{i}\in {\cal{B}}'(K)\}$. Next, we will construct the nerve of the closed coverings ${\cal{B}}_{K}$ and
${\cal{B}}'(K)$. For this purpose we will always consider open balls even if the collections consist of closed balls.
Let $N_{{\cal{B}}'(K)}$ and $N_{{\cal{B}}_{K}}$ be the corresponding geometric realization of the nerves of ${\cal{B}}'(K)$ and ${\cal{B}}_{K}$, respectively.
For simplicity, we will call $N_{{\cal{B}}'(K)}$ and $N_{{\cal{B}}_{K}}$ {\it the nerves} of ${\cal{B}}'(K)$ and ${\cal{B}}_{K}$, respectively.\\

\begin{lem}
Let $N_{{\cal{B}}'(K)}$, $N_{{\cal{B}}_{K}}$ be the nerves of ${\cal{B}}'(K)$ and ${\cal{B}}_{K}$,
respectively. Then  $N_{{\cal{B}}'(K)}$ is homeomorphic to $N_{{\cal{B}}_{K}}$  and homeomorphic to $K$.
\end{lem}

{\it Proof.} Recall from section 3.1 that the nerve in ${\mathbb R}^{n+2}$ of balls in ${\cal{B}}(K)$ is a collection of $(n+2)$-simplexes and $(n+1)$-dimensional
diamonds with empty interior. By adding a ball $B_{C}$ at a center $C$ of a flower (which is in turn the center of a diamond) the effect on the nerve is to add
the cone from $C$, thus filling the interior of the diamond. This is the nerve of the collection ${\cal{B}}'(K)$.\\

Let $S_{K}(v)$, $S_{{\cal{B}}'(K)}(v)$ and $S_{{\cal{B}}_{K}}(v)$ denote the star of a vertex $v$ in the simplicial complexes  $K$, $N_{{\cal{B}}'(K)}$ and
 $N_{{\cal{B}}_{K}}$ respectively. We will define an isotopy between these simplicial complexes locally on $S_{K}(v)$, leaving the boundary
$Lk_{K}(v)$ of $S_{K}(v)$ fixed.\\

Let $P$ be a $n$-hyperplane supporting the $n$-dimensional face of a cube $E\in {\cal{C}}_{2}$. Let $v$ be a vertex in $K\cap P$. If  $S_{K}(v)$ is contained
in $P$, the centers of all the balls surrounding $v$ lie on $P$, and therefore  $S_{K}(v)$, $S_{{\cal{B}}'(K)}(v)$ and $S_{{\cal{B}}_{K}}(v)$ all coincide and
are equal to the $n$-disk which is the geometric convex hull of the vertices surrounding $v$ (see Figure 5).\\

\begin{figure}[tbh]
\centerline{\epsfxsize=1.6in \epsfbox{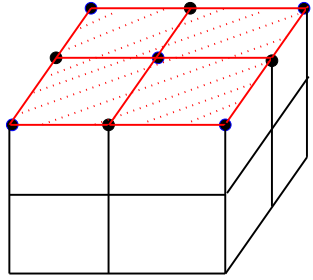}}
\caption{\sl A schematic picture of the corresponding geometric nerve of balls centered at $F$.}
\label{F5}
\end{figure}

Now take a vertex $v$ which is a center of a flower on an $(n-1)$-dimensional cubic simplex $S$ in which $K$ turns. Then $S_{K}(v)$ is an
$n$-dimensional disk whose boundary contains vertices $v_{1},\ldots,v_{r}$ on different faces of a cube $E\in {\cal{C}}_{2}$. Observe that if
$v_{1},\ldots,v_{r}$ belong to more than two faces of $E$ then $v$ has to be a vertex of the cube and, inversely, all the vertices of $E$ are centers
of flowers. By construction we did not add a ball with center $v$ when we constructed ${\cal{B}}'(K)$. By proposition 4.3, only the balls forming a flower in $v$
intersect. Therefore $N_{{\cal{B}}'(K)}$ is, locally around $v$, formed by $n$-dimensional faces of a diamond with vertices in $v_{1},\ldots,v_{r}$, and these
faces contain the boundary of  $S_{K}(v)$ (see Figure 6).\\

\begin{figure}[tbh]
\centerline{\epsfxsize=1.6in \epsfbox{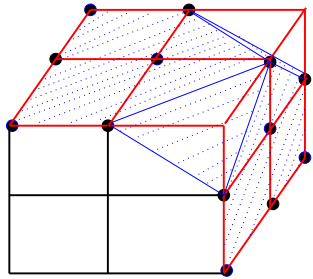}}
\caption{\sl A schematic picture of the corresponding geometric nerve of balls centered at $F_{1}\cup F_{2}$.}
\label{F6}
\end{figure}

We can isotope $S_{K}(v)$ into $N_{{\cal{B}}'(K)}$ leaving this boundary fixed. After this isotopy, all
of $K$, $N_{{\cal{B}}'(K)}$ and $N_{{\cal{B}}_{K}}$ coincide. $\blacksquare$\\

We will use the following theorem to prove that $\tilde{T}$ is a regular neighborhood of $K$.

\begin{theo} Let ${\mathcal{B}}:=\{B_1,
\cdots,{B_r}\}$ be a finite set of round open balls in either ${\Bbb R}^{n+2}$ with its euclidean metric or ${\Bbb S}^{n+2}$
with its spherical metric. Let $v_i$ denote the center of $B_i$. Let $\bar{B_i}$ ($i=1\cdots{r}$) denote the corresponding
closed balls and $\Bbb S^{n+1}_i:=\partial\bar{B}_i$ their spherical boundaries. We will assume that if two balls are tangent then
there exists $j$ such that $B_j$ is centered at the point of tangency and that none of balls is contained in another.  Suppose that
the {\it geometric nerve} of ${\mathcal{B}}$ consisting of totally geodesic simplexes (with respect to the euclidean or spherical metric,
 respectively) defined before (see section 3.1) is an $n$-dimensional polyhedral sphere $K$ in ${\Bbb R}^{n+2}$ or ${\Bbb S}^{n+2}$
such that for each closed ball $\bar{B_i}$
the pair $(\bar{B_i}, \bar{B_i}\cap{K})$ is equivalent to the trivial
tangle (in particular $K$ is locally flat). Then ${\mathcal{N}}(K):=\cup_{i=i}^{r}\,\,\bar{B_i}$ is a closed regular
neighborhood of $K$ and, in particular ${\mathcal{N}}(K)\cong{{\Bbb D}^2}\times{K}$ since the normal bundle of an $n$-dimensional
locally flat sphere in ${\Bbb R}^{n+2}$ or ${\Bbb S}^{n+2}$ is trivial.

\end{theo}

{\it Proof.}  Let $V(K)$ be a regular neighborhood of $K$. Since $K$ is locally flat,
there exists a homeomorphism $\eta:{{\Bbb D}^2}\times{K}\to{V(K)}$. For $0<\delta\leq1$, let
$V_\delta(K):=\eta({\Bbb D^2}_\delta\times{K})$, where ${\Bbb D^2}_\delta$ is the closed 2-disk of radius $\delta$.\\

The proof will continue by induction on the dimension $m=n+2$. When $m=2$,  $\mathcal{B}$ consists of two balls $B_1$ and $B_2$
(in $\Bbb R^2$ or $\Bbb S^2$) such that $\bar{B_1}\cap\bar{B_2}=\emptyset$ and the result is obviously true. If $m=3$ we have a pearl
necklace in $\Bbb R^3$ or $\Bbb S^3$ \cite{homogeneity} and the proof is contained there (see Figure 7).

\begin{figure}[tbh]
\centerline{\epsfxsize=1.25in \epsfbox{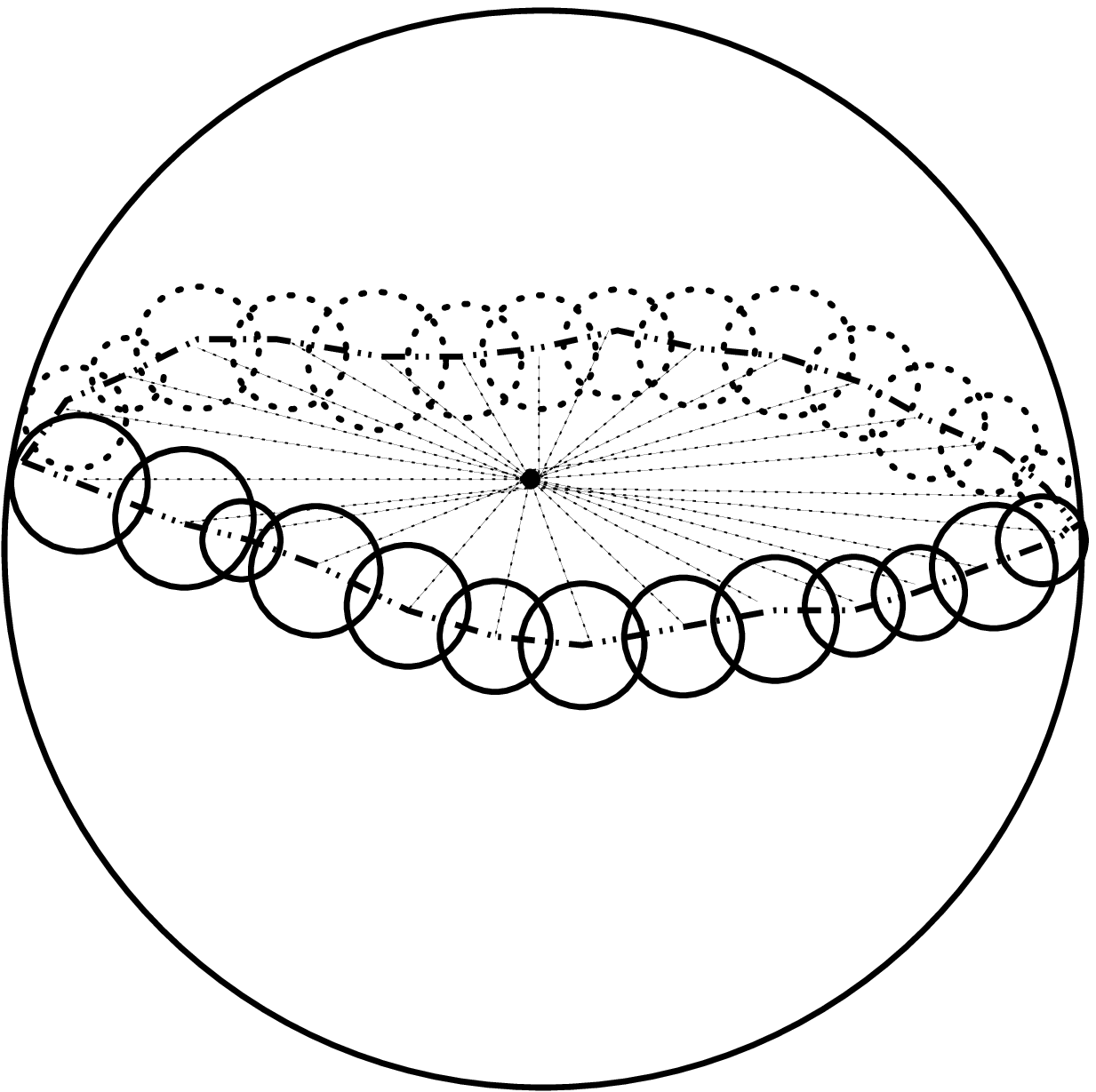}}
\caption{\sl A schematic pearl necklace of dimension 3.}
\label{F7}
\end{figure}

Let $m>3$. Let  $j\in\{1,\ldots,r\}$, and ${\mathcal{B}}_j:=\{B_i\,\,\,|\,B_i\cap{B_j}\neq\emptyset\}$.  The geometric realization of
the nerve of ${\mathcal{B}}_j$, which we denote by $S_j$,
is the star of the vertex corresponding to the center $v_j$ of the ball $B_j$.
The set of spherical $(n+1)$-balls $L_j:=\{\bar{B}_k\cap{\Bbb S^{n+1}_j}\,\,|\,\,\bar{B_k}\cap{\bar{B}_j}\neq\emptyset\}$
satisfies the induction hypotheses of the theorem applied to the ambient space $\Bbb S^{n+1}_j$ and, in fact, the geometric realization
of these $(n+1)$-balls in $\Bbb S^{n+1}_j$ is combinatorial equivalent ({\it i.e.,} PL-homeomorphic) to the link of the vertex $v_j$ 
{\it i.e., }the boundary of the star $S_j$.\\

By hypothesis, the pair $(\bar{B}_j,\bar{B}_j\cap{K})$ is equivalent to the trivial $n$-dimensional tangle.
Let $M_j:=\cup_{\bar{B}_k\in{L_k}}\,\,\bar{B}_k\cap{\bar{B_j}}$.
For a fixed $j\in\{1,\ldots,r\}$ let $V_j:=M_j\cup\,\,(V_\delta(K)\cap\bar{B_j})$ (see Figure 8).
Then for $\delta$ sufficiently small $V_j$ is a tubular neighborhood of $\bar{B}_j\cap{K}$ in $\bar{B}_j$. More precisely
the pair $(\bar{B}_j,V_j)$ is homeomorphic as a pair to $(I^{n+2},I^n\times\frac12{I^2})$ where $I^{n+2}=[-1,1]\times\cdots\times[-1,1]$
(($n+2$)-times) and $I^n\times\frac12{I^2}=[-1,1]\times\cdots\times[-1,1]\times[-\frac12,\frac12]\times[\frac12,\frac12]$ (the factor $[-1,1]$
occurs $n$-times). Since,  $\bar{B}_j-\text{Int}\,(V_j)\cong{I^{n+2}-\text{Int}\,(I^n)}\times\frac12{I^2}\cong\partial{I^2}\times[-1,1]\times{I^n}$
one has that $V_j$ is a strong deformation retract of $\bar{B}_j$.\\

\begin{figure}[tbh]
\centerline{\epsfxsize=1.25in \epsfbox{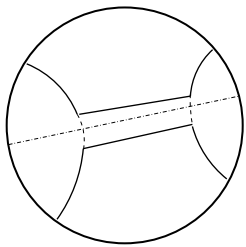}}
\caption{\sl A schematic picture of the tubular neighborhood $V_j$.}
\label{F8}
\end{figure}

The previous arguments imply that the set $\cup_{j=1}^{r}{V_j}$ is a regular neighborhood of $K$ and furthermore
since $\cup_{j=1}^{r}{V_j}$ is a strong deformation retract
of ${\mathcal{N}}(K)=\cup_{j=1}^{r}\bar{B}_j$ it follows that  ${\mathcal{N}}(K)$ is also a regular neighborhood of $K$. $\blacksquare$\\

\begin{coro}
Let $\tilde{T}$ be a generalized increased pearl necklace subordinate to $K$. Then $\tilde{T}$ is isotopically equivalent to a closed regular neighborhood $\cal{N}$ of $K$.
\end{coro}

\subsection{Description of the limit set}

Let $K$ be a smooth $n$-knot. Let $T=\cup^{r}_{i=1}B_{i}$ be a generalized pearl necklace subordinate to $K$ and $\tilde{T}$ be the corresponding generalized increased
pearl necklace.
Let $\Gamma$ be the group generated by reflections $I_{j}$ through
$\partial B_{j}$,
$B_{j}\in T$. To guarantee that
the group $\Gamma$ is Kleinian we will use the Poincar\'e Polyhedron Theorem (see \cite{epstein}, \cite{kap2}, \cite{maskit}).
This theorem establishes conditions for the group to be discrete. In practice
these conditions are very hard to verify, but in our case all of them are
satisfied automatically from the construction, since the balls $B_{i}$, $B_{j}\in T$ are either disjoint, tangent or
their boundaries meet orthogonally.\\

This theorem also gives us a presentation for the group $\Gamma$. In our case, the dihedral angles $n_{ij}$
between the faces $F_{i}$, $F_{j}$ are $\frac{\pi}{2}$, if the faces are adjacent or $0$ otherwise (by definition). Therefore,
we have the following presentation of $\Gamma$
$$
\Gamma=<I_{j}, j=1,\dots,n|\hspace{.2cm}(I_{j})^{2}=1,\hspace{.2cm}(I_{i}I_{j})^{n_{ij}}=1>
$$

The fundamental domain for $\Gamma$ is $D=\mathbb{S}^{n+2}-T$ and $\Gamma$ is geometrically finite.\\

The natural question for the Kleinian group $\Gamma$ is: What is its limit set?
Recall that
to find the limit set of $\Gamma$, we need to find all the accumulation points of its orbits. We will do this by stages.\\

\hspace{-.67cm}{\bf Stage I.} We will apply induction on the number of reflections.
\begin{enumerate}
\item {\it First step}: We reflect with respect to $\partial B_{1}$, {\it i.e.}, we apply $I_{1}\in\Gamma$. Notice that
$D_{1}=D\cup I_{1}(D)=(\mathbb{S}^{n+2}-T)\cup (B_{1}-(I_{1}(T))$ is a fundamental domain of an index two
subgroup $\Gamma_{1}$ in $\Gamma$. \\

{\it Claim:} The set $D^{'}_{1}=(\tilde{T}-B_{1})\cup  I_{1}(\tilde{T})$ is a
regular neighborhood of the geometric realization of its nerve, which is isotopic  to the
connected sum of $K$ with its mirror image $-K$.

In fact, observe that the reflection map $I_{1}$ replaces the trivial tangle
$(B_{1}, B_{1}\cap K)$ by a new tangle  $C_{1}=(B_{C_{1}},K_{C_{1}})$
 which is isotopic to an
orientation-reversing copy of the tangle 
$$
C=(\mathbb{S}^{n+2}-\mbox{Int}(B_{1}), K-\mbox{Int}(B_{1}\cap K))=(B_{C},K_{C}),
$$
which is homeomorphic to the canonical tangle associated to $K$ up to isotopy. Then we apply theorem 4.5. 

\item {\it Second step}: We reflect with respect to $\partial B_{2}$, {\it i.e.}, we apply $I_{2}\in\Gamma$. Notice that
$D_{2}=D_{1}\cup I_{2}(D_{1})=(\mathbb{S}^{n+2}-T)\cup (B_{1}-I_{1}(T))\cup (B_{2}-I_{2}(D_{1}))$ is a
fundamental domain of an index two
subgroup $\Gamma_{2}$ in $\Gamma_{1}$.

The set
$D^{'}_{2}=(\tilde{T}-\cup _{i=1}^{2}B_{i})\cup I_{1}(\tilde{T})\cup I_{2}(D^{'}_{1})$ is a
regular neighborhood of the geometric realization of its nerve, which is isotopic to the
connected sum $K\#(-K)\#K\#(-K)$.

\item {\it $r^{th}$-step}:  We reflect with respect to $\partial B_{r}$, {\it i.e.}, we apply $I_{r}\in\Gamma$. Notice that
$D_{r}=D_{r-1}\cup I_{r}(D_{r-1})$ is a
fundamental domain of an index two
subgroup $\Gamma_{r}$ in $\Gamma_{r-1}$.

The set
$D^{'}_{r}=I_{1}(\tilde{T})\cup I_{2}(D_{1}^{'})\cup I_{3}(D_{2}^{'})\cup\ldots\cup I_{r}(D_{r-1}^{'})$ is a
regular neighborhood of the geometric realization of its nerve, which is isotopic  to the
connected sum of $2^{r-1}$ copies of $K$ and $2^{r-1}$ copies of $-K$.

\end{enumerate}

At the end of the $r^{th}$-step , we obtain a regular neighborhood $\tilde{T}_{1}$ of a new tame knot $K_{1}$, which is isotopic
to the connected sum of $2^{r-1}$ copies of $K$ and $2^{r-1}$ copies of $-K$. Notice that
$\tilde{T}_{1}\subset\mbox{Int}(\tilde{T})$.\\

\hspace{-.67cm}{\bf Stage II.}\\

Repeat $k$-times Stage I. At this stage we obtain a regular neighborhood $\tilde{T}_{k}$ of a new tame knot $K_{k}$ which is isotopic to
the connected sum of $2^{kr-1}$ copies of $K$ and $2^{kr-1}$ copies of $-K$. By construction,
$\tilde{T}_{k}\subset \mbox{Int}(\tilde{T}_{k-1})$.\\

Let $x\in\cap_{l=1}^{\infty} \tilde{T}_{l}$. We shall prove that $x$ is a
 limit point. Indeed, there exists a sequence of closed balls
 $\{B_{m}\}$ with $B_{m}\subset \tilde{T}_{m}$ such that $x\in B_{m}$ for
 each $m$. We can find a $z\in\mathbb{S}^{n+2}-{\tilde{T}}$
 and a sequence $\{w_{m}\}$ of distinct elements of $\Gamma$, such
 that $w_{m}(z)\in B_{m}$. Since $diam(B_{m})\rightarrow 0$ it follows
 that $w_{m}(z)$ converges to $x$. Hence $\cap_{l=1}^{\infty} {\tilde{T}}_{l}\subset \Lambda (\Gamma)$. The other inclusion
$\Lambda (\Gamma)\subset \cap_{l=1}^{\infty} \tilde{T}_{l}$ clearly holds.
Therefore, the limit set is given by
$$
\Lambda (\Gamma)=\varprojlim_{l} \tilde{T}_{l}=\bigcap_{l=1}^{\infty} \tilde{T}_{l}.
$$
\vskip .3cm
\begin{rem}
This description is similar to the one which appears in the work of Peter Scott on \cite{S}: Suppose that $G$ is a right-angled reflection
group with the fundamental domain $P$ in $\mathbb{H}^{n}$. Given a facet $F$ of $P$, let $P_{F}$ denote the union
of $P$ and its reflection in $F$. Then $P_{F}$ is a fundamental domain of an index two subgroup in $G$. Now, define inductively
index two subgroups $G=G_{0}>G_{1}>G_{2}>\ldots$, where $G_{i}$ has the fundamental domain $P_{i}=PF_{-1}$. Then
$\cap_{i}G_{i}=\{1\}$. In particular, if we apply this to our construction, the union of the fundamental domains is
the entire discontinuity set.
\end{rem}

\begin{theo}
Let $T$ be a generalized pearl necklace
 of a  non-trivial tame knot $K$ of dimension $n$ consisting of closed round balls in $\cal{B}$. Let $\tilde{T}$ be the
corresponding generalized increased pearl necklace.
 Let $\Gamma$ be the group generated by reflections on the boundary of each ball of $T$.
Let $\Lambda(\Gamma)$ be the corresponding limit set. Then
 $\Lambda(\Gamma)$ is homeomorphic to $\mathbb{S}^{n}$.
\end{theo}

\hspace{-.67cm}{\it Proof.}
The proof makes use of an infinite process similar to the one used to prove that the Whitehead manifold cross $\mathbb{R}$ is
homeomorphic to $\mathbb{R}^{4}$ (\cite{mcmillan}, \cite{poenaru}).\\

Each regular neighborhood $\tilde{T}_{i}$ of the $n$-dimensional knot $K_{i}$ is homeomorphic to
$K_i\times\mathbb{D}^2$ and satisfies
$\tilde{T}_{i}\subset \tilde{T}_{i-1}$, for
each $i$. Note that $K_{i}$ is knotted in $\tilde{T}_{i-1}$, in the sense that $K_{i-1}\times{0}$ is not isotopic in $\tilde{T}_{i-1}$ to
$K_{i}$. \\

Consider $\tilde{T}_{i-1}\times I$, $I=[0,1]$. Since $K_{i}$ is of codimension 3 in  $\tilde{T}_{i-1}\times I$ then it is not
{\it topologically} knotted in
$\tilde{T}_{i-1}\times I$ (\cite{brakes}, \cite{stallings}), {\it i.e.}, it is isotopic, in $\tilde{T}_{i-1}\times I$, to $K_{i-1}\times0$. \\

The fact that $K_i$ is not topologically knotted in $\tilde{T}_{i-1}\times I$ implies that $\tilde{T}_{i}\times I$ can be deformed
 isotopically inside  $\tilde{T}_{i-1}\times I$ onto a small tubular neighborhood of $K_{i-1}$.\\

The above facts imply that the pair $\tilde{T}_{m+1}\times \frac{1}{2^{m+1}}I\subset \tilde{T}_{m}\times \frac{1}{2^{m}}I$ is
topologically equivalent to
$\mathbb{S}^{n}\times\frac{1}{2^{m+1}}\mathbb{D}^{3}\subset \mathbb{S}^{n}\times\frac{1}{2^{m}}\mathbb{D}^{3}$, {\it i.e.}, we have
the following commutative diagram
$$
\begin{CD}
\tilde{T}_{1}\times I@<<<\tilde{T}_{2}\times \frac{1}{2}I@<<<\cdots \tilde{T}_{m}\times \frac{1}{2^{m}}I@<<<\cdots\Lambda\\
@V\sim VV@V\sim VV@V\sim VV@V VV\\
\mathbb{S}^{n}\times\mathbb{D}^{3}@<<<\mathbb{S}^{n}\times\frac{1}{2}\mathbb{D}^{3}@<<<\cdots\mathbb{S}^{n}\times\frac{1}{2^{m}}\mathbb{D}^{3}
@<<<\cdots\mathbb{S}^{n}\\
\end{CD}
$$

Therefore, by the universal property of the inverse limits, there exists a homeomorphism from $\Lambda$ to $\mathbb{S}^{n}$ associated
to this sequence of maps. $\blacksquare$\\

\section{Dynamically-defined fibered wild knots}

 In 1925 Emil Artin described two methods for constructing knotted spheres of
dimension $n$ in $\mathbb{S}^{n+2}$ from knots in $\mathbb{S}^{n+1}$. One of these methods is called {\it spinning}
and  uses the rotation process.
A way to visualize it is the following. We can consider $\mathbb{S}^{2}$ as an
$\mathbb{S}^{1}$-family of half-equators (meridians=$\mathbb{D}^{1}$) such that
the respective points of their boundaries are identified to obtain the poles. Then the formula
$Spin(\mathbb{D}^{1})=\mathbb{S}^{2}$ means to send
homeomorphically the unit interval $\mathbb{D}^{1}$ to a meridian of
$\mathbb{S}^{2}$ such that $\partial\mathbb{S}^{1}=\{0,1\}$ is mapped to the poles
and, multiply the interior of $\mathbb{D}^{1}$ by $\mathbb{S}^{1}$. In other words, one
spins the meridian with respect to the poles to obtain $\mathbb{S}^{2}$.
Similarly, consider $\mathbb{S}^{n+1}$ as an $\mathbb{S}^{1}$-family of half-equators
($\mathbb{D}^{n}$) where boundaries are respectively identified, hence
$Spin(\mathbb{D}^{n})=\mathbb{S}^{n+1}$ means to send homeomorphically
$\mathbb{D}^{n}$ to a meridian of $\mathbb{S}^{n+1}$ and keeping $\partial\mathbb{D}^{n}$
fixed, multiply the interior of $\mathbb{D}^{n}$ by $\mathbb{S}^{1}$. In particular, if we start with a fibered
tame knot $K$, then $Spin(K)$ also fibers (see \cite{zeeman1}).\\

We recall that a knot  $K$ in $\mathbb{S}^{n+2}$ is {\it fibered} if there exists a
locally trivial fibration $f:(\mathbb{S}^{n+2}-K)\rightarrow \mathbb{S}^{1}$. We
require further that $f$ be well-behaved near $K$, that is, that it
has a neighborhood framed as $\mathbb{S}^{n}\times\mathbb{D}^{2}$, with
$K\cong \mathbb{S}^{n}\times\{0\}$, in such a way that the
restriction of $f$ to $\mathbb{S}^{n}\times(\mathbb{D}^{2}-\{0\})$ is the map
into $\mathbb{S}^{1}$ given by $(x,y)\rightarrow \frac{y}{|y|}$. It follows that each
$f^{-1}(x)\cup K$, $x\in\mathbb{S}^{1}$, is a ($n+1$)-manifold
 with boundary $K$: in fact a Seifert (hyper-) surface for $K$
(see \cite{rolfsen}, page 323).

\begin{exs}
\begin{enumerate}
\item The right-handed trefoil knot and the figure-eight knot are fibered knots with fiber
the punctured torus.
\item The trivial knot $\mathbb{S}^{n}\hookrightarrow\mathbb{S}^{n+2}$ is fibered by the projection map
$(\mathbb{S}^{n}*\mathbb{S}^{1})-\mathbb{S}^{n}\rightarrow\mathbb{S}^{1}$. Fibers are (n+1)-disks.
\item Let $K\subset\mathbb{S}^{3}$ be a fibered tame knot. Then $Spin(K)$ also fibers. Using again the Spin process,
we can get fibered knots in any dimension. Thus, there exist non-trivial fibered knots in any dimension $n\geq3$ (compare \cite{ansu}).
\end{enumerate}
\end{exs}

Next, we will apply  our construction to a non-trivial fibered tame knot $K$ of dimension $n$, with fiber  $S$.
 Let $T=\cup^{r}_{i=1}B_{i}$ be a generalized pearl necklace subordinate to $K$ and $\tilde{T}$ be the corresponding generalized increased
pearl necklace. Let $\Gamma$ be the group generated by reflections
$I_{j}$ on the boundary
of balls $B_{j}\in T$, $j=1,\ldots, r$. Observe  that the fundamental domain for the group $\Gamma$ is
$D=\mathbb{S}^{n+2}-\mbox{Int}(T)$ (see definition 2.4) and in this case $D$ is homeomorphic to $\Omega (\Gamma)/\Gamma$
(see section 2).

\begin{lem}
 Let $T=\cup^{r}_{i=1}B_{i}$ be a generalized pearl necklace subordinate to
 a  non-trivial fibered tame knot $K$ of dimension $n$, with fiber  $S$. Then
$\Omega (\Gamma)/\Gamma$ fibers over the circle with fiber $S^{*}$,
diffeomorphic to the closure of $S$ in $\mathbb{S}^{n+2}$.
\end{lem}

\hspace{-.67cm}{\bf {\it Proof.}}
Let $\widetilde{P}:(\mathbb{S}^{n+2}-K)\rightarrow \mathbb{S}^{1}$ be the
given fibration with fiber the manifold $S$. Observe that
$\widetilde{P}\mid_{\mathbb{S}^{n+2}-\mbox{Int}(T)}\overset{\rm def}= P$ is a fibration and, after
modifying $\widetilde{P}$ by isotopy if necessary, we can consider  that
the fiber $S$  cuts the boundary of each ball $B_{i}\in T$ transversely, in $n$-disks (see Figure 9).\\

\begin{figure}[tbh]
\centerline{\epsfxsize=1.2in \epsfbox{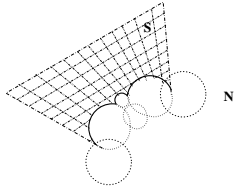}}
\caption{\sl A schematic picture of the fiber intersecting each ball in an $n$-disk.}
\label{F10}
\end{figure}

Hence the space $D$ fibers over the
circle with fiber the $(n+1)$-manifold $S^{*}$, which is the closure of the
manifold $S$ in $\mathbb{S}^{n+2}$ (see section 2). Since
$\Omega(\Gamma)/\Gamma\cong D$, the
result follows. $\blacksquare$\\

By the above Lemma, in order to describe completely the orbit space
$(\mathbb{S}^{n+2}-\Lambda(\Gamma))/\Gamma$
in the case that the original knot
is fibered, we only need to determine its monodromy, which is
precisely the monodromy of the knot.\\

Consider the orientation-preserving index two subgroup $\widetilde{\Gamma}\subset\Gamma$.
The fundamental domain for $\widetilde{\Gamma}$ is
$$
\widetilde{D}=D\cup I_{j}(D)=(\mathbb{S}^{n+2}-\mbox{Int}(T))\cup (B_{j}-\mbox{Int}(I_{j}(T))),
$$
for some $I_{j}\in\Gamma$ and $B_{j}\in T$. Observe that $I_{j}(D)\cap D\cong \mathbb{S}^{n-1}\times\mathbb{D}^{2}$.

\begin{lem}
Let $T=\cup^{r}_{i=1}B_{i}$ be a generalized pearl necklace subordinate to
 a  non-trivial fibered tame knot $K$ of dimension $n$, with fiber  $S$. Then
$\Omega (\widetilde{\Gamma})/\widetilde{\Gamma}$ fibers over the circle with fiber a
$(n+1)$-manifold $S^{**}$, which is homeomorphic to the ($n+1$)-manifold
$S^{*}$ joined along an n-disk to a copy of itself in $\mathbb{S}^{n+2}$ modulo $\widetilde{\Gamma}$.
\end{lem}

\hspace{-.67cm}{\bf {\it Proof.}} First, we will prove that the fundamental domain $\widetilde{D}$ for $\widetilde{\Gamma}$ fibers over
the circle.

Let $\widetilde{P}:(\mathbb{S}^{n+2}-K)\rightarrow \mathbb{S}^{1}$ be the
given fibration with fiber the manifold $S$ and
let $P:(\mathbb{S}^{n+2}-T)\rightarrow \mathbb{S}^{1}$ be its corresponding restriction.
Observe that the canonical tangle $K_{T}$ associated to $K$ also fibers
over the circle with fiber $S$, hence $B_{j}-I_{j}(K)$ fibers over the circle with fiber $S_{j}$ which is homeomorphic to $S$.
 By the same argument of the above Lemma, $I_{j}(D)=B_{j}-I_{j}(T)$ fibers over the circle
with fiber $S_{j}^{*}$,
which is homeomorphic to $S^{*}$, the closure of $S$, via the fibration $P_{j}$.\\

Given $\theta\in\mathbb{S}^{1}$, let $P^{-1}(\theta)=S^{*}_{\theta}$ and $P_{j}^{-1}(\theta)=S_{j\theta}^{*}$ be the corresponding
fibers. Notice that $I_{j}(S^{*}_{\theta})=S_{j\theta}^{*}$, hence
$S^{*}_{\theta}\cap\partial B_{j}=S_{j\theta}^{*}\cap\partial B_{j}$. Therefore $\widetilde{D}$ fibers over the circle with fiber
a $(n+1)$-manifold $\hat{S}_{\theta}$, which is homeomorphic to the ($n+1$)-manifold
$S^{*}$ joined along an $n$-disk to a copy of itself in $\mathbb{S}^{n+2}$.\\

Since $\Omega (\widetilde{\Gamma})/\widetilde{\Gamma}$ is homeomorphic to $\widetilde{D}/\widetilde{\Gamma}$ (see section 2), and
$\widetilde{D}/\widetilde{\Gamma}$ fibers over the circle with fiber $S^{**}=\hat{S}/\sim$, where
$\partial B_{K}\sim I_{j}(\partial B_{k})$ $k\neq j$, via $I_{k}I_{j}^{-1}\in\widetilde{\Gamma}$. The result follows.
$\blacksquare$\\

Since
$\widetilde{\Gamma}$ is a normal subgroup of $\Gamma$, it follows by Lemma 8.1.3 in
\cite{thurston1} that $\widetilde{\Gamma}$ has the same limit set than $\Gamma$. Therefore
$\mathbb{S}^{n+2}-\Lambda(\Gamma)=\mathbb{S}^{n+2}-\Lambda(\widetilde{\Gamma})$.

\begin{theo}
Let $T=\cup^{r}_{i=1}B_{i}$ be a generalized pearl necklace subordinate to
a  non-trivial fibered tame knot $K$ of dimension $n$, with fiber $S$ . Let $\Gamma$ be the group generated
by reflections on each ball of $T$ and let $\widetilde{\Gamma}$
be the orientation-preserving index two subgroup of $\Gamma$.
Let $\Lambda(\Gamma)=\Lambda(\widetilde{\Gamma})$ be the
corresponding limit set. Then:
\begin{enumerate}
\item There exists a locally trivial fibration $\psi
:(\mathbb{S}^{n+2}-\Lambda(\Gamma))\rightarrow\mathbb{S}^{1}$, where the
 fiber $\Sigma_{\theta}=\psi^{-1}(\theta)$ is an orientable ($n+1$)-manifold with one end, which is homeomorphic to the connected
sum along $n$-disks of an infinite number of copies of $S$.
\item  $\overline{\Sigma_{\theta}}-\Sigma_{\theta}=\Lambda(\Gamma)$, where
$\overline{\Sigma_{\theta}}$ is the closure of $\Sigma_{\theta}$ in $\mathbb{S}^{n+2}$.
\end{enumerate}
\end{theo}

\hspace{-.67cm}{\bf {\it Proof.}}
We will first prove that $\mathbb{S}^{n+2}-\Lambda(\Gamma)$ fibers over the circle.
We know that $\zeta:\Omega(\widetilde{\Gamma)}\rightarrow\Omega(\widetilde{\Gamma)}/\widetilde{\Gamma}$
is an infinite-fold covering since $\widetilde{\Gamma}$ acts freely on $\Omega(\widetilde{\Gamma})$. By
the previous lemma, there exists a locally trivial
fibration $\phi:\Omega(\widetilde{\Gamma})/\widetilde{\Gamma}\rightarrow\mathbb{S}^{1}$ with fiber
$S^{**}$.
Then $\psi=\phi\circ\zeta:\Omega(\widetilde{\Gamma})\rightarrow\mathbb{S}^{1}$ is a
locally trivial fibration. The fiber is $\Gamma(S^{*})$, {\it i.e.}, the orbit of
the fiber.  \\

We now describe $\Sigma_{\theta}=\Gamma(S^{*})$ in detail.
Let $\widetilde{P}:(\mathbb{S}^{n+2}-K)\rightarrow \mathbb{S}^{1}$ be the
given fibration. The fibration $\widetilde{P}\mid_{\mathbb{S}^{n+2}-\mbox{Int}(T)}\overset{\rm def}= P$
has been chosen as in the lemmas above. The
fiber $\widetilde{P}^{-1}(\theta)=P^{-1}(\theta)$ is a Seifert surface $S^{*}$ of $K$, for each
$\theta\in\mathbb{S}^{1}$. We suppose $S^{*}$ is oriented. Recall that the boundary of $S^{*}$ cuts
each ball $B_{j}$ in an $n$-disk $a_{j}$.\\

The reflection  $I_{j}$ maps both a copy of $T-B_{j}$
(called $T^{j}$) and a copy of $S^{*}$ (called $S_{1}^{*j}$) into
the ball $B_{j}$, for $j=1,2,\ldots,r$. Observe that both $T^{j}$ and $S_{1}^{*j}$ have opposite orientation
and that $S^{*}$ and $S_{1}^{*j}$ are joined by
the $n$-disk $a_{j}$ (see Figure 10) which, in both manifolds, has the same orientation.\\

\begin{figure}[tbh]
\centerline{\epsfxsize=.8in \epsfbox{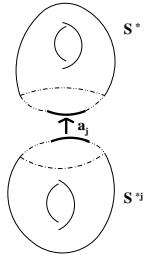}}
\caption{\sl A schematic picture of the sum of two Seifert surfaces $S^{*}$ and $S^{*j}$ along the $n$-disk $a_{j}$.}
\label{F11}
\end{figure}

Recall the description of the limit set $\Lambda$ of $\Gamma$ in section 4. At the end of the first stage, we
have a new regular neighborhood $\tilde{T}_{1}$ of the knot $K_{1}$. By the same argument as above, its complement fibers
over the circle with fiber the Seifert surface $S^{*}_{1}$, which is in turn
homeomorphic to the sum of $N_{1}$ copies of $S^{*}$ along the respective $n$-disks.\\

Continuing this process, at the end we obtain that $\Lambda(\Gamma)$ fibers over the circle with fiber
$\Sigma_{\theta}$ which is homeomorphic to the connected sum along $n$-disks of an infinite number of copies of $S$. Notice that
the diameter of the $n$-disks tends to zero.\\

Next, we will describe its set of ends.
Consider the Fuchsian model (see \cite{maskit}). In this case, we are considering the trivial knot. Then its limit set
is the unknotted sphere $\mathbb{S}^{n}$ and its complement fibers over $\mathbb{S}^{1}$
with fiber the disk $\mathbb{D}^{n+1}$.
\vskip .3cm
In each step we are adding copies of $S$ to this disk in such a way that
they accumulate on the boundary. If we intersect this disk with any
compact set, we have just one connected component. Hence it has only
one end. Therefore, our Seifert surface has one end (see Figure 11).

\begin{figure}[tbh]
\centerline{\epsfxsize=1in \epsfbox{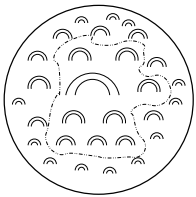}}
\caption{\sl A schematic picture of a disk with copies of $S$ intersected with a compact set.}
\label{F12}
\end{figure}

The first part of the theorem has been proved. For the second part,
observe that the closure of the fiber in $\mathbb{S}^{n+2}$ is the fiber
union its boundary. Therefore
$\overline{\Sigma_{\theta}}-\Sigma_{\theta}=\Lambda(\Gamma)$. $\blacksquare$

\begin{rem}
This theorem gives an open book decomposition (\cite{winkelnkempe}, \cite{rolfsen}
pages  340-341) of
$\mathbb{S}^{n+2}-\Lambda(\Gamma)$, where the ``binding'' is the
knot $\Lambda(\Gamma)$, and each ``page'' is an
orientable ($n+1$)-manifold with one end (the fiber).\\
Indeed, this decomposition can be thought of in the following way. By the above
theorem, $\mathbb{S}^{n+2}-\Lambda(\Gamma)$ is $\Sigma_{\theta}\times [0,1]$ modulo the
identification of the top with the bottom through a characteristic
homeomorphism. Consider $\overline{\Sigma_{\theta}}\times [0,1]$ and identify
the top with the bottom. This is equivalent to keeping
$\partial\overline{\Sigma_{\theta}}$ fixed
and spinning $\Sigma_{\theta}\times \{0\}$ with respect to
$\partial\overline{\Sigma_{\theta}}$ until it is
glued  with $\Sigma_{\theta}\times \{1\}$.
Removing $\partial\overline{\Sigma_{\theta}}$ we obtain the open book decomposition.
\end{rem}

\section{Monodromy}

Let $K$ be a non-trivial tame fibered $n$-knot and let $S$ be the fiber. Since
$\mathbb{S}^{n+2}-K$ fibers over the circle, we know that
$\mathbb{S}^{n+2}-K$ is a mapping torus equal to $S\times [0,1]$ modulo
a characteristic homeomorphism $\psi:S\rightarrow S$ that glues $S\times\{0\}$ to
$S\times\{1\}$. This homeomorphism induces a homomorphism
$$
\psi_{\#}:\Pi_{1}(S)\rightarrow\Pi_{1}(S)
$$
called {\it the monodromy of the fibration}.\\

Another way to understand the monodromy is through  the
{\it  Poincar\'e's first return map} of a flow, defined as follows. Let $M$ be
connected, compact
manifold and let $f_{t}$ be a flow that possesses a transverse
section $\eta$. It follows that if $x\in\eta$ then there exists a
continuous function $t(x)>0$ such
that $f_{t}\in\eta$. We may define Poincar\'e's first return map
$F:\eta\rightarrow\eta$ as $F(x)=f_{t(x)}(x)$. This map is a diffeomorphism
and induces a homomorphism of $\Pi_{1}$
called {\it the monodromy} (see \cite{verjovsky}, chapter 5).\\

For the manifold $\mathbb{S}^{n+2}-K$, the flow that defines the Poincar\'e's first
return map $\Phi$ is the flow that cuts transversely each page of
its open book decomposition.\\

Consider a generalized pearl necklace $T$ subordinate to $K$. As we have observed
during the reflecting process, $K$ and $S$ are copied in each
reflection. So the flow $\Phi$ is also copied. Hence, the
Poincar\'e's map can be extended at each stage, giving us in the end a
homeomorphism $\psi:\Sigma_{\theta}\rightarrow \Sigma_{\theta}$ that identifies
$\Sigma_{\theta}\times\{0\}$ with $\Sigma_{\theta}\times\{1\}$. This homeomorphism  induces
the monodromy of the limit $n$-knot.\\

By the long exact sequence associated to a fibration, we have
\begin{equation}
0\rightarrow \Pi_{1}(\Sigma_{\theta})\rightarrow
\Pi_{1}(\mathbb{S}^{n+2}-\Lambda (\Gamma ))\stackrel{\Psi}{\overset{\longleftarrow}{\rightarrow}}
\mathbb{Z}\rightarrow 0, \tag{1}
\end{equation}
which has a homomorphism section
$\Psi:\mathbb{Z}\rightarrow\Pi_{1}(\mathbb{S}^{n+2}-\Lambda (\Gamma))$.
Therefore (1) splits.
As a consequence
$\Pi_{1}(\mathbb{S}^{n+2}-\Lambda (\Gamma))$ is
the semi-direct product of $\mathbb{Z}$ with $\Pi_{1}(\Sigma_{\theta})$.
This gives a method for computing the fundamental group of a limit $n$-knot
whose complement fibers over the circle.\\

Notice that there is only one homomorphism from
$\Pi_{1}(\mathbb{S}^{n+2}-\Lambda (\Gamma))$ onto $\mathbb{Z}$, since by Alexander duality $H^{1}(\mathbb{S}^{n+2}-\Lambda (\Gamma))\cong\mathbb{Z}$.
Therefore, the monodromy of the limit knot $\Lambda (\Gamma)$ is completely determined by the monodromy of the knot $K$.

\begin{rem}
The monodromy provides us with a way to distinguish between two limit fibered knots $\Lambda_{1}$ and $\Lambda_{2}$.
\end{rem}

\begin{theo} Let $T$ be a generalized pearl necklace
 of a  non-trivial tame fibered $n$-knot $K=Spin^{n-1}(K')$, where $K'\subset\mathbb{S}^{3}$ is a
non-trivial tame fibered knot of dimension one and $n=1,\ldots,5$.
Let $\Gamma$ be the group generated by reflections on each ball of $T$.
Let $\Lambda(\Gamma)$ be the corresponding limit set. Then
 $\Lambda(\Gamma)$ is wildly embedded in
 $\mathbb{S}^{n+2}$.
\end{theo}

{\it Proof}. Let $K'$ be a non-trivial tame fibered 1-knot. Then the
fiber $S$ is a Seifert surface of genus $g$ whose boundary is $K'$.
The
fundamental group of $S$ is the free group in $2g$ generators $\{a_{i},\,b_{i}\,:\,i=1,\ldots,g\}$, hence
the fundamental group of the fiber $\Sigma_{\theta}$ of the limit $n$-knot is the free group  with generators  $\{a^{j}_{i},\,b^{j}_{i}\,:\,i=1,\ldots,g;\,\,j\in\mathbb{N}\}$.
Since $\Pi_{1}(\mathbb{S}^{n+2}-K)\cong\Pi_{1}(\mathbb{S}^{3}-K')$ (see \cite{rolfsen}), it follows that the
monodromy maps, in both cases, coincide. Let $\psi_{\#}$ be the monodromy of $K$. Then
the monodromy of the limit $n$-knot
$\widetilde{\psi_{\#}}:\Pi_{1}(\Sigma_{\theta})\rightarrow
\Pi_{1}(\Sigma_{\theta})$ sends $a^{j}_{i}\mapsto (\psi_{\#}(a_{i}))^{j}$ and $b^{j}_{i}\mapsto (\psi_{\#}(b_{i}))^{j}$.\\

Thus,
$$
\begin{aligned}
\Pi_{1}(\mathbb{S}^{n+2}-\Lambda (\Gamma))&\cong
\Pi_{1}(\mathbb{S}^{1})\ltimes_{\widetilde{\psi_{\#}}}\Pi_{1}(\Sigma_{\theta})\\
&=\{a^{j}_{i},\,b^{j}_{i},\,c\,:\,a^{j}_{i}*c=(\psi_{\#}(a_{i}))^{j},\hspace{.2cm}b^{j}_{i}*c=(\psi_{\#}(b_{i}))^{j}\}.\\
\end{aligned}
$$

Therefore, the fundamental group $\Pi_{1}(\mathbb{S}^{n+2}-\Lambda (\Gamma))$ is infinitely generated.
This implies that $\Lambda (\Gamma)$ is wildly embedded. $\blacksquare$

\begin{coro}
There exist infinitely many nonequivalent wild $n$-knots in $\mathbb{R}^{n+2}$.
\end{coro}

\begin{coro}Let $T$ be a generalized pearl-necklace whose template is a
non-trivial tame fibered $n$-knot $K$. Then
$\Pi_{1}(\Omega(\Gamma)/\Gamma)\cong\mathbb{Z}\ltimes_{\psi_{\#}}\Pi_{1}(\Sigma_{\theta})$.
\end{coro}

As a consequence of Theorem 4.8, Theorem 6.2 and Corollary 6.3, we have:\\

\hspace{-.67cm}{\bf THEOREM 4.1}
{\it There exist infinitely many non-equivalent knots $\psi:\mathbb{S}^{n}\rightarrow \mathbb{S}^{n+2}$ wildly embedded
as limit sets of geometrically finite Kleinian groups, for $n=1,2,3,4,5$.}\\

\noindent{\bf Acknowledgments.} We would like to thank Cynthia Verjovsky \hyphenation{Marcotte} Marcotte
for her suggestions after carefully reading our paper. We would also like to thank the referee for his/her valuable suggestions.

M. Boege. Instituto de Matem\'aticas, Unidad Cuernavaca. Universidad Nacional Au\-t\'o\-no\-ma de M\'exico.
Av. Universidad s/n, Col. Lomas de Chamilpa. Cuernavaca, Morelos, M\'exico, 62209.

\hspace{-.6cm}{\it E-mail address:} margaret@matcuer.unam.mx
\vskip .3cm
G. Hinojosa. Facultad de Ciencias, Universidad Aut\'onoma del Estado de Morelos. Av. Universidad 1001, Col. Chamilpa.
Cuernavaca, Morelos, M\'exico, 62209.

\hspace{-.6cm}{\it E-mail address:} gabriela@buzon.uaem.mx
\vskip .3cm
A. Verjovsky. Instituto de Matem\'aticas, Unidad Cuernavaca. Universidad Nacional Au\-t\'o\-no\-ma de M\'exico.
Av. Universidad s/n, Col. Lomas de Chamilpa. Cuernavaca, Morelos, M\'exico, 62209.

\hspace{-.6cm}{\it E-mail address:} alberto@matcuer.unam.mx

\end{document}